\documentclass[smallextended]{svjour3}       
\usepackage{graphicx}
\usepackage{amssymb}
\usepackage{amsmath}
\begin{document}

\title{A Family of Elliptic Curves With Rank $\geq5$
}


\author{Farzali Izadi         \and
        Kamran Nabardi 
}


\institute{Farzali Izadi \at
              Department of Mathematics, Azarbaijan Shahid Madani University, Tabriz 53751-71379, Iran \\
              \email{farzali.izadi@azaruniv.edu}           
           \and
           Kamran Nabardi \at
              Department of Mathematics, Azarbaijan Shahid Madani University, Tabriz 53751-71379, Iran\\
              \email{nabardi@azaruniv.edu}
}
\maketitle
\begin{abstract}
In this paper, we construct a family of elliptic curves with rank $\geq 5$. To do this, we use the  Heron formula for a triple $(A^2, B^2, C^2)$ which are not necessarily the three sides of a triangle. It turns out that  as parameters of a family of elliptic curves, these three positive integers $A$, $B$, and $C$, along with the extra parameter $D$ satisfy the quartic Diophantine equation $A^4+D^4=2(B^4+D^4)$.
\end{abstract}

\section{Introduction}
\label{intro}
As is well-known, the affine part of an elliptic curve $E$ over a field $\Bbb{K}$ can
be explicitly expressed by the generalized Weierstrass equation of the form
\begin{equation}
E:  y^2+a_1xy+a_3y=x^3+a_2x^2+a_4x+a_6,\end{equation}
where $a_1,a_2,a_3,a_4,a_6\in \Bbb{K}$.  In this paper we are interested
in the case of
$\Bbb{K}=\Bbb{Q}$.\\
By the Mordell-Weil theorem \cite{Was}, every elliptic curve over
$\Bbb{Q}$ has a commutative group $E(\Bbb{Q})$ which is finitely
generated, i.e., $E(\Bbb{Q})\cong\Bbb{Z}^r\times
E(\Bbb{Q})_{tors} $, where $r$ is a nonnegative integer called the rank of $E(\Bbb{Q})$ and
$E(\Bbb{Q})_{tors}$ is the subgroup of elements of finite order
called the torsion subgroup of $E(\Bbb{Q})$.\\
By the Mazur theorem \cite{Sil1}, the torsion subgroup $E(\Bbb{Q})_{tors}$ is
one of the following 15 types: $\mathbb{Z}/n\mathbb{Z}$ with
$1\leq n\leq 10$ or $n=12$, $\mathbb{Z}/2\mathbb{Z}\times
\mathbb{Z}/2m\mathbb{Z}$ with $1\leq m\leq 4$. Besides, it is not
known which values of rank $r$ are possible. The folklore
conjecture is that a rank can be arbitrarily  large, but it seems to
be very difficult to find examples with large ranks. The current record
is an example of elliptic curve over $\mathbb{Q}$ with rank $\geq$
28, found by Elkies  in May 2006 (see \cite{Duj}). Having
classified the torsion part, one interested in seeing whether or
not the rank is unbounded among all the elliptic curves. There is
no known guaranteed algorithm to determine the rank and it is not
known which numbers can occur as the ranks.
\section{Previous Works}\label{sec2}
Let $D$ be a non-zero integer. The curve
\begin{equation*}
E_D\ :\ y^2=x^3+Dx
\end{equation*}
has been considered by many mathematicians. Note that the  congruent number elliptic curve belongs to this category by  $D=-n^2$. By taking $D=p$, where $p\equiv5\pmod{8}$ and less than $1000$, Bremner and Cassels \cite{Bremner} show  that the rank is always 1 in accordance with the conjecture of Selmer and Mordell. Kudo and Motose \cite{Kudo} studied the case $D=-p$, where $p$ is a Fermat or Mersenne prime and found the ranks $0$, $1$ and $2$. By taking $D=pq$, where $p$ and $q$ are distinct odd primes, Yoshida \cite{Yoshida} showed that the rank is at most $5$. In \cite{IKN}, the authors take $D=-n$, where $n=u^4+v^4=r^4+s^4$ and prove the rank is at least $3$. Moreover, they show that if $n$ is odd and the parity conjecture is true, then the rank is at least $4$. Maenishi \cite{Maenishi} studied the case $D=-pq$, where $p$ and $q$ are distinct odd prime numbers and by imposing  an extra condition found a rank $4$ family.\\

 In  the case of $D=-n^2$ ( congruent number elliptic curve),  many authors have attempted to find curves with high ranks. Rogers \cite{Rogers1}, using an idea of Rubin and Silverberg \cite{Rubin}, found two curves with ranks $5$ and $6$. Later, he found \cite{Rogers2} a curve with rank $7$. In \cite{Duje1}, by using Mestre-Nagao's sum \cite{Duje2,Nagao1,Nagao2}, the authors wish  to find  congruent number elliptic curves with high ranks. They succeed to find new congruent elliptic curves with rank $6$. Johnstone \cite{Jennifer}, in her Master thesis provides an in depth background on congruent numbers and elliptic curves and then, presents a family of congruent number elliptic curves with rank at least three.\\

In this work  we consider the elliptic curve
\begin{equation*}
E:\ y^2=x^3-4S^2x,
\end{equation*}
over the $K3$ surface
\begin{equation*}
T=V_{1,1,-2-2}:\ A^4+D^4=2(B^4+C^4),
\end{equation*}
where, $S=S(A,\ B,\ C,\ D)$ is a rational function of $A$, $B$, $C$ and $D$. We prove that the group of rational maps
$P:\ T\longrightarrow E$,
that commute with the projection $E\longrightarrow T$, has rank at least $5$. We do this by exhibiting five explicit sections $P_1$, $P_2$, $P_3$, $P_4$, and $P_5$, and later showing that they are independent. We use the fact that $T(\mathbb{Q})$ is infinite to deduce that infinitely many  specializations of $E$ have the rank at least $5$ over $\mathbb{Q}$. This is done by using some elliptic curves of positive rank lying on $T$  that we found in \cite{IN,IN1}.

 A natural question is that whether the set of rational points $T(\mathbb{Q})$ is Zariski dense in the $K3$ surface $T$. The following theorem answers the question affirmatively.
 \begin{theorem}
 Let $a, b, c, d\in\mathbb{Q}^\ast$ be nonzero rational numbers with $abcd$ square. Let $P=[x_0\ :\ y_0\ :\ z_0\ :\ w_0]$ be a rational point on $V_{a, b, c, d}:\ ax^4+by^4+cz^4+dw^4=0$, and suppose that $x_0y_0z_0w_0\not=0$ and let $P$ is not contained in one of the $48$ lines of the surface. Then the set of rational points of $V$ is Zariski dense in $V$ as well as dense  in  the real analytic topology on $V(\mathbb{R})$ .
 \end{theorem}
 \begin{proof}
 \cite{Luijk}, Theorem 1.1.
 \end{proof}

 The point $(21, 19, 20, 7)$ is on $V_{1,1,-2,-2}$ and satisfies the hypotheses of above theorem, so it shows that the rational points on the surface $V_{1,1,-2-2}$ are indeed dense in Zariski topology,  and also in the real analytic topology.
 \section{Definition of rational function $S$}\label{sec3}

In this work we deal with a family of elliptic curves which are related to  the positive
integer solutions of the diophantine equation
\begin{equation}\label{Eq1.0}
A^4+D^4=2(B^4+C^4)
\end{equation}
as follows.\\
 Heron formula states that for a triangle with sides $a$, $b$ and $c$, one can get the area of the triangle by the formulae:
\begin{equation}\label{Eq1.1}
S=\sqrt{p(p-a)(p-b)(p-c)},
\end{equation}
where, $p=\frac{a+b+c}{2}$.\\

   Take  $(A^2, B^2, C^2)$, where   $(A, B, C)$  are as in  equation \eqref{Eq1.0}. Since the triple $(A, B, C)$ is arising from the Diophantine equation \eqref{Eq1.0}, there is not guarantee to have a real triangle. Now by taking
       $a=A^2$, $b=B^2$, and $c=C^2$, we find that
\begin{equation}\label{Eq22}
S=\sqrt{\frac{(A^2+B^2+C^2)(A^2+B^2-C^2)(A^2+C^2-B^2)(B^2+C^2-A^2)}{16}}.
\end{equation}
Expanding \eqref{Eq22}, one gets
\begin{equation*}
S^2=-\left(\frac{A^8+B^8+C^8-2A^4B^4-2A^4C^4-2B^4C^4}{16}\right),
\end{equation*}
equivalently,
\begin{equation*}
16S^2=2A^4B^4+2A^4C^4+2B^4C^4-A^8-B^8-C^8,
\end{equation*}
or
\begin{equation*}
\left(\frac{A^4+B^4-C^4}{2}\right)^2+4S^2=A^4B^4.
\end{equation*}

\begin{equation*}
\left(\frac{A^4+B^4-C^4}{2}\right)^2+4S^2=A^4B^4.
\end{equation*}
Multiplying both sides by $A^2B^2$ yields
\begin{equation*}
A^2B^2\left(\frac{A^4+B^4-C^4}{2}\right)^2+4A^2B^2S^2=A^6B^6.
\end{equation*}
Taking $y=AB\left(\frac{A^4+B^4-C^4}{2}\right)$ and $x=A^2B^2$, we get the following  family of elliptic curves:
\begin{equation}\label{Eq1.2}
E:\ y^2=x^3-4S^2x.
\end{equation}

Since the roles of $A$, $B$, and $C$ are symmetric in the Heron formula, we have the following points on
the family \eqref{Eq1.2} also.
\begin{equation*}
\begin{array}{l}
P_1=\left(A^2B^2,\ \frac{AB(A^4+B^4-C^4)}{2}\right),\\
P_2=\left(A^2C^2,\ \frac{AC(A^4+C^4-B^4)}{2}\right),\\
P_3=\left(B^2C^2,\ \frac{BC(B^4+C^4-A^4)}{2}\right),
\end{array}
\end{equation*}
Next we wish to impose two more  points on the curve \eqref{Eq1.2} with $x-$coordiates as $B^2D^2$ and $C^2D^2$. Substituting $x=B^2D^2$ in \eqref{Eq1.2} yields
\begin{equation*}
y^2=B^2D^2\left(\frac{4B^4D^4+A^8+B^8+C^8-2A^4B^4-2A^4C^4-2B^4C^4}{4}\right),
\end{equation*}
or
\begin{equation}\label{Eq1.3}
y^2=B^2D^2\left(\frac{A^4(A^4-2B^4-2C^4)+B^8+C^8-2B^4C^4+4B^4D^4}{4}\right).
\end{equation}
Let
\begin{equation}\label{Eq1.4}
A^4-2B^4-2C^4=-D^4,
\end{equation}
then by substituting $A^4=2B^4+2C^4-D^4$ in \eqref{Eq1.3}, we get
\begin{equation*}
y^2=B^2D^2\left(\frac{B^2+D^2-C^2}{2}\right)^2.
\end{equation*}

Thus, the point $P_4=\left(B^2D^2,\ \frac{BD(B^4+D^4-C^4)}{2}\right)$ is a new  point on \eqref{Eq1.2}.
Similarly, one can easily check that  the point $P_5=\left(C^2D^2,\ \frac{CD(C^4+D^4-B^4)}{2}\right)$ lies  on \eqref{Eq1.2} as well.\\

It is clear that the existence  of these extra points on the family depends exactly on the existence of the solutions of the Diophantine equation $A^4+D^4=2(B^4+C^4)$.\\

In the next stage, we will talk about the torsion subgroup of the \eqref{Eq1.2}.

\section { The Torsion subgroup of \eqref{Eq1.2}}\label{sec4}
Before starting the argument on  the torsion subgroup of \eqref{Eq1.2}, let us recall the following theorem which reveals the structure of torsion subgroup for elliptic curves of the form $y^2=x^3+Dx$, where $D\in\Bbb{Z}$ and is fourth- power-free integer.
\begin{theorem}\label{t2}
Let $D\in\Bbb{Z}$ be a fourth-power-free integer, and $E_D$ be the elliptic curve
\begin{equation*}
E_D:\ y^2=x^3+Dx.
\end{equation*}
Then,
\begin{equation*}
E_{D,tors}(\Bbb{Q})\simeq\left\{\begin{array}{lll}
\Bbb{Z}/4\Bbb{Z},&&\text{if}\ D=4,\\
&&\\
\Bbb{Z}/2\Bbb{Z}\times\Bbb{Z}/2\Bbb{Z},&& \text{if}\ -D\  \text{is a perfect square},\\
&&\\
\Bbb{Z}/2\Bbb{Z},&&\text{otherwise}.
\end{array}\right.
\end{equation*}
\end{theorem}
\begin{proof}
See \cite[Proposition 6.1, page 346]{Sil}
\end{proof}

\begin{remark}\label{remark1}
 From the appearance of \eqref{Eq1.2} and  also Theorem \ref{t2}, one may guess that the  family defined in \eqref{Eq1.2} is a congruent number elliptic curve family and so, the torsion group is $\Bbb{Z}/2\Bbb{Z}\times\Bbb{Z}/2\Bbb{Z}$. In the following we show this is not the case.
 First of all, one should be aware of the probable values of $S$ in \eqref{Eq22}.  As $(A^2, B^2, C^2)$  does not  necessarily  construct a triangular, there is no guarantee for  $S$ to be a real or an imaginary number in equation \eqref{Eq22}.  For example, for  the 4-tuple $(A,D,B,C)=(21,19,20,7)$, we get $S=180\sqrt{979}$, while $(A,D,B,C)=(1661081,988521,336280,1437599)$ leads us to
 \begin{equation*}
 S=840\sqrt{962357334498800500956761065836542898196489}\ i.
 \end{equation*}
  In the former case $S$ is real, while in the later one $S$ is an imaginary number. In the   case of   real valued  $S$, we see that  the coefficient of $x$ in \eqref{Eq1.2} is negative. But  if $S$ is an imaginary number, the coefficient of  $x$  would be positive in \eqref{Eq1.2}. Moreover, $S=180\sqrt{979}$ shows that it is not necessary  for $4S^2$  to be a perfect square. In this case, whether or not the value of $4S^2$ be a perfect square, the torsion subgroup of  \eqref{Eq1.2} is $\Bbb{Z}/2\Bbb{Z}\times\Bbb{Z}/2\Bbb{Z}$ or $\Bbb{Z}/2\Bbb{Z}$.
 If $S$ is an imaginary number, then the coefficient of $x$ in \eqref{Eq1.2} is positive. Therefore, the torsion subgroup of \eqref{Eq1.2} is $\Bbb{Z}/2\Bbb{Z}$.
\end{remark}

\section{Finding the solutions of  the equation $A^4+D^4=2(B^4+C^4)$}\label{sec5}
In \cite{IN} Izadi and Nabardi found infinitely many integer solutions of this equation. Their method is based on the points of the elliptic curve $y^2=x^3-36x$  with a generator $(-3, \ 9)$  explained as follows:\\
 Let  $P_n=(x_n,\ y_n)$, where $P_n=n\cdot(-3,\ 9)$ ($n\in\mathbb{N}$) is a point on the  elliptic curve $y^2=x^3-36x$, one gets
\begin{equation}\label{Eq1.5}
\begin{array}{l}
A_n={\phi_n}^4 + 1296{{\psi}_n}^8 + 864{\phi_n}{{\psi}_n}^6 + 72{\phi_n}^2{{\psi}_n}^4 + 144{\omega_n}{{\psi}_n}^5 \\
  \qquad\qquad - 24{\phi_n}^3{{\psi}_n}^2 + 4{\phi_n}^2{\omega_n}{{\psi}_n},\\
  \\
 D_n=-864{\phi_n}{{\psi}_n}^6 - {\phi_n}^4 - 1296{{\psi}_n}^8 - 72{\phi_n}^2{{\psi}_n}^4 + 144{\omega_n}{{\psi}_n}^5\\
  \qquad\qquad + 24{\phi_n}^3{{\psi}_n}^2 + 4{\phi_n}^2{\omega_n}{{\psi}_n},\\
  \\
  B_n=4({\phi_n}^2 + 36{{\psi}_n}^4){\omega_n}{{\psi}_n},\\
  \\
  C_n=({\phi_n}^2 - 36{{\psi}_n}^4 - 12{\phi_n}{{\psi}_n}^2)({\phi_n}^2 - 36{{\psi}_n}^4 + 12{\phi_n}{{\psi}_n}^2),\\
\end{array}
\end{equation}
 such that $A_n^4+D_n^4=2(B_n^4+C_n^4)$ ,  where $\phi_n$ and $\psi_n$ are $n$-th division polynomials ( see \cite[pp.80-83]{Was} ). Therefore, we define a family of elliptic curves by
\begin{equation}\label{Eq1.6}
E_n:\ y^2=x^3+\left(\frac{A_n^8+B_n^8+C_n^8-2A_n^4B_n^4-2A_n^4C_n^4-2B_n^4C_n^4}{4}\right)x.
\end{equation}
As we discussed above, there are $5$ points on \eqref{Eq1.6} which are given by:
\begin{equation*}
\begin{array}{l}
P_{n1}=\left(A_n^2B_n^2,\ \frac{A_nB_n(A_n^4+B_n^4-C_n^4)}{2}\right),\\
\\
P_{n2}=\left(A_n^2C_n^2,\ \frac{A_nC_n(A_n^4+C_n-B_n^4)}{2}\right),\\
\\
P_{n3}=\left(B_n^2C_n^2,\ \frac{B_nC_n(B_n^4+C_n^4-A_n^4)}{2}\right),\\
\\
P_{n4}=\left(B_n^2D_n^2,\ \frac{B_nD_n(B_n^4+D_n^4-C_n^4)}{2}\right),\\
\\
P_{n5}=\left(C_n^2D_n^2,\ \frac{C_nD_n(C_n^4+D_n^4-B_n^4)}{2}\right).\\
\end{array}
\end{equation*}

Let $n=2$, then
\begin{equation*}
E_2:\ y^2=x^3+2716157340889414533900362432217058675869770553600x,
\end{equation*}
and
\begin{equation*}
\begin{array}{l}
P_{21}=\left(110502951275524201934400,\ 549083548316905650689533416877852800\right),\\
\\
P_{22}=\left(2019516118036966895564241,\ -3704296107487960167032542005050395239\right),\\
\\
P_{23}=\left(23710164715943220558400,\ 804710464588380886496762950328163200\right),\\
\\
P_{24}=\left( 312020909765749236942400,\ 936950008965894699667383086290460800\right),\\
\\
P_{25}=\left(5702393005462282638861361,\ 141744467546800549687092185696272575\right).
\end{array}
\end{equation*}
Using SAGE  \cite{Sage} software, shows that the determinant of height pairing matrix of $[P_{21},P_{22},\ P_{23},\ P_{24},\ P_{25}]$ equals $\cong30739535.349$ and so, these points are independent. Since specialization is an injective homomorphism \cite[pp.456-457]{Sil}, it follows that the equation \eqref{Eq1.6} is a family of elliptic curves with rank $\geq 5$.\\

In a separate paper \cite{IN1}, the authors  found another set of solutions for  \eqref{Eq1.4}. Let us  briefly recall  this method.
The smallest known solution for \eqref{Eq1.4} is  $(A_0,D_0,B_0,C_0)=(21,19,20,7)$. Dividing  the equation \eqref{Eq1.4} by $C^4$, we get
\begin{equation}\label{EQ5}
x^4+y^4-2u^4-2=0\qquad x,y,u\in\Bbb{Q}.
\end{equation}
Now $(x_0,y_0,u_0)=(3,19/7,20/7)$ is a solution for \eqref{EQ5}. In this stage we define
\begin{equation}\label{EQ6}
x=at+x_0,\qquad y=bt+y_0,\qquad u=ct+u_0.
\end{equation}
Putting these in \eqref{EQ5}, yields
\begin{equation}\label{EQ7}
Mt^4+Nt^3+Rt^2+St=0,
\end{equation}
where
\begin{equation}\label{EQ8}
\begin{array}{l}
M=a^4+b^4-2c^4,\\
\\
N=12a^3+\frac{76}{7}b^3-\frac{160}{7}c^3,\\
\\
R=54a^2+\frac{2166}{49}b^2-\frac{4800}{49}c^2,\\
\\
S=108a+\frac{27436}{343}b-\frac{64000}{343}c.
\end{array}
\end{equation}
By letting $S=0$, one has
\begin{equation}\label{EQ9}
c=\frac{9261}{16000}a+\frac{6859}{16000}b.
\end{equation}
If we let $R=0$, then  from \eqref{EQ9} we have
\begin{equation}\label{EQ10}
b=\frac{12147}{10507}a,\ \frac{93}{133}a.
\end{equation}
Let $b=\frac{93}{133}a$, so we can take $c=\frac{123}{140}a$. From \eqref{EQ7} we have $t=-\frac{N}{M}$ and so
\begin{equation}\label{EQ11}
t=-\frac{1732800}{389209a}.
\end{equation}
Therefore the new rational solution for \eqref{EQ5} is:
\begin{equation}\label{EQ12}
x_1=-\frac{565173}{389209},\qquad y_1=-\frac{1086621}{2724463},\qquad u_1=-\frac{2872540}{2724463}.
\end{equation}
Equivalently,
\begin{equation}\label{EQ13}
A_1=-3956211,\qquad D_1=-1086629,\qquad B_1=-2872540,\qquad C_1=2724463
\end{equation}
is a new solution for \eqref{Eq1.4}. Again by repeating the same argument, we can find another solution for \eqref{Eq1.4}. Corresponding to $(607, 1999, 951, 1640)$, we have
\begin{equation}
y^2=x^3+9749352988442901002400000x,
\end{equation}
using MWRANK \cite{Mwrank}, one can show the rank is $8$.\\
Let $(181, 2077, 1247, 1620)$, then we get
\begin{equation}
y^2=x^3+4988940634912192616750400x,
\end{equation}
In this case, the rank is $6$.\\




\end{document}